\documentclass[a4paper,reqno]{amsart}

\usepackage{amssymb}
\usepackage{amsmath}
\usepackage{tikz-cd}
\usepackage{hyperref}

\newtheorem{theorem}{Theorem}
\newtheorem{corollary}{Corollary}

\theoremstyle{definition}
\newtheorem{example}{Example}
\newtheorem{remark}{Remark}

\newcommand{\cZ}{\mathcal{Z}}

\DeclareMathOperator{\Ext}{Ext}
\DeclareMathOperator{\Hilb}{Hilb}
\DeclareMathOperator{\Pic}{Pic}
\DeclareMathOperator{\Hom}{Hom}
\DeclareMathOperator\IZ{\mathcal{I}_{\mathcal{Z}}}
\DeclareMathOperator\Db{D^b}
\newcommand\Sn{S^{[n]}}

\title{Smooth components on special iterated Hilbert schemes}

\author{Fabian Reede}
\address{Institut f\"ur Algebraische Geometrie, Leibniz Universit\"at Hannover, Welfengarten 1, 30167 Hannover, Germany}
\email{reede@math.uni-hannover.de}



\begin{document}

\begin{abstract}
	Let $S$ be a smooth projective surface with $p_g=q=0$. We show how to use derived categorical methods to study the geometry of certain special iterated Hilbert schemes associated to $S$ by showing that they contain a smooth connected component isomorphic to $S$.
\end{abstract}

\maketitle

\vspace{-2ex}
\section{Introduction}
Hilbert schemes are ubiquitous in modern algebraic geometry. But even in good situations these schemes can behave badly. This became clear with Mumford's famous example, which shows that there is an irreducible component of the Hilbert scheme of smooth irreducible curves in $\mathbb{P}^3$ of degree 14 and genus 24 that is generically non-reduced, see \cite{mum}. More exactly Mumford constructs a 56-dimensional irreducible family of such curves, such that the tangent space at each point of this component has dimension 57 and proves that this family is not contained in any other irreducible family of dimension bigger than $56$.

In this note we prove that there are certain iterated Hilbert schemes which contain at least one smooth connected component. More exactly the main result of this note is:
\begin{theorem}\label{thm}
	Assume $S$ is a smooth projective surface with $p_g=q=0$ and let $\Sn$ be the Hilbert scheme of length $n$ subschemes of $S$. Then the universal family $\cZ$ in $S\times \Sn$ can be understood as a family of codimension two subschemes in $\Sn$ with common Hilbert polynomial $p(t)$ classified by $S$ such that the classifying morphism identifies $S$ with a smooth connected component of the Hilbert scheme $\Hilb^{p(t)}(\Sn)$.
\end{theorem}
This theorem has its origin in a result by Lange and Newstead, who proved a similar result for curves and moduli spaces of stable vector bundles on these curves in \cite{lang}. More exactly they show that if $M$ is a fine moduli space of stable vector bundles on a smooth projective curve $C$ of genus $g\geqslant 2$ with universal family $\mathcal{U}$ on $C\times M$, then for any $c\in C$ the vector bundle $\mathcal{U}_c$ on $M$ is stable. Furthermore they show that for $c\neq c'$ we have $\mathcal{U}_c\not\cong\mathcal{U}_{c'}$. Together with a previous result by Narasimhan and Ramanan these results imply that $C$ embeds as a smooth connected component in a moduli spaces of stable vector bundles o $M$.

The main input into the proof of the main result of this note is a result by Krug and Sosna which states that the integral functor $\Phi: \Db(S) \rightarrow \Db(\Sn)$ with kernel the universal ideal sheaf $\IZ$ is fully faithful. This result allows to reduce the computation of certain $\Ext$-groups on $\Sn$ to the computation of easier $\Ext$-groups on $S$.

All objects in this note are defined over the field of complex numbers $\mathbb{C}$.

\subsection*{Acknowledgement} I thank Pieter Belmans for informing me about the fully faithfulness results in \cite{bel}  and \cite{krug} as well as Ziyu Zhang for many useful conversations.

\section{Proof of the Main Theorem}

Let $S$ be a smooth projective surface with $p_g=q=0$, that is we have
\begin{equation*}
\mathrm{H}^1(S,\mathcal{O}_S)=\mathrm{H}^2(S,\mathcal{O}_S)=0.
\end{equation*}
\begin{remark}
Around the year 1870 Max Noether posed the question if surfaces with $p_g=q=0$ are necessarily rational, see \cite[\S 3, Question 1]{bauer}. By now it is known that the answer to this question is negative. In fact besides rational surfaces there is a huge class of surfaces satisfying these conditions which are not rational, most classically Enriques surfaces which have Kodaira dimension zero. But there are also surfaces of general type satisfying these conditions for example Godeaux surfaces, Campedelli surfaces or Beauville surfaces, see \cite{bauer2,dolga} for more information and examples.
\end{remark}

In the following we denote the Hilbert scheme of length $n$ subschemes of $S$ by $\Sn$, that is we have as sets:
\begin{equation*}
	\Sn=\left\lbrace [Z]\,\,|\,\, Z\subset S\,\, \text{is a zero-dimensional subscheme with}\,\, h^0(Z,\mathcal{O}_Z)=n \right\rbrace.
\end{equation*}
It is well known that $\Sn$ is smooth and that $\dim(\Sn)=2n$. Using this notation we have the universal subscheme
\begin{equation}\label{univ}
	\cZ=\left\lbrace (s,[Z])\in S\times \Sn\,\,|\,\, s\in\text{supp}(Z) \right\rbrace \subset S\times \Sn 
\end{equation}
coming with the corresponding universal ideal sheaf $\IZ\hookrightarrow \mathcal{O}_{S\times\Sn}$.

\begin{remark}\label{flat}
Recall that the universal family $\cZ$ is flat over $\Sn$. Indeed, using definition \eqref{univ} one can see that the restriction of $p: S\times \Sn \rightarrow \Sn$ to $\cZ$ is finite and flat of degree $n$. But as a matter of fact $\cZ$ is also flat over $S$ due to \cite[Theorem 2.1]{krug1}.
\end{remark}

As we use integral functors in the following, we quickly recall their definition: let $X$ and $Y$ be smooth projective varieties and denote their bounded derived categories of coherent sheaves by $\Db(X)$ and $\Db(Y)$ respectively, then the integral functor with kernel $\mathcal{K}\in \Db(X\times Y)$ is defined by
\begin{equation*}
	\Phi_\mathcal{K}: \Db(X)\rightarrow \Db(Y),\,\,\, E\mapsto \mathbf{R}p_{*}(q^{*}E\otimes^L \mathcal{K})
\end{equation*}
where $p$ and $q$ are the projections $X\times Y \rightarrow Y$ resp. $X\times Y\rightarrow X$, see \cite[\S 5]{huy2}.

The description of the image of an integral transform for a skyscraper sheaf $\mathcal{O}_x$ of a closed point $x\in X$ is rather easy, which we will collect as
\begin{example}\cite[Examples 5.4 (vi)]{huy2}\label{exam}
	Assume the kernel of the integral functor $\Phi_{\mathcal{K}}$ is in fact a coherent sheaf on $X\times Y$ flat over $X$, then we have for every closed point $x\in X$
	\begin{equation*}
		\Phi_\mathcal{K}(\mathcal{O}_x)\cong \mathcal{K}_x,
	\end{equation*}
where the fiber $\mathcal{K}_x:=\mathcal{K}_{|\left\lbrace x\right\rbrace \times Y}$ is considered as a sheaf on $Y$ via the second projection $\left\lbrace x\right\rbrace \times Y \rightarrow Y$.
\end{example}

Interpreting the universal ideal sheaf $\IZ$ as an element in $\Db(S\times\Sn)$ (as a complex concentrated in degree zero), we can look at the integral functor with kernel given by $\IZ$:
\begin{equation*}
\Phi: \Db(S) \rightarrow \Db(\Sn),\,\,\,E \mapsto \mathbf{R}p_{*}(q^{*}E\otimes \IZ).
\end{equation*}
 In our case this integral functor has some very good properties. The main input into this note is the following very useful fact discovered by Krug and Sosna:
\begin{theorem}\cite[Theorem 1.2]{krug}\label{faith}
	Let $S$ be a smooth projective surface which satisfies $p_g=q=0$, then the integral functor $\Phi$ is fully faithful.
\end{theorem}

As an application of Theorem \ref{faith} we have the following
\begin{corollary}\label{ext}
	Assume $S$ is a smooth projective surface with $p_g=q=0$ then for all $E,F\in\Db(S)$ and $i\geqslant0$ there is an isomorphism
	\begin{equation*}
		\Ext^i_{\Sn}(\Phi(E),\Phi(F))\cong \Ext^i_S(E,F).
	\end{equation*}
\end{corollary}

\begin{remark}
There are also fully faithfulness results for universal families of moduli spaces of stable bundles of rank two and degree one on smooth projective curves of genus $g\geqslant 2$ by Narasimhan as well as Fonarev and Kuznetsov, see \cite{nara1}, \cite{nara2} and \cite{fon}. Recently these results were generalized to higher rank and degree by Belmans and Mukhopadhyay as well as Lee and Moon, see \cite{bel3} and \cite{moon}. These results can be used to give a proof of the result of Lange and Newstead in the spirit of this note.
\end{remark}

\textbf{Proof of Theorem}\,\ref{thm}\textbf{.} Remark \ref{flat} shows that the universal family $\cZ\subset S\times \Sn$ is flat over $S$. Denote the fiber over a closed point $s\in S$ by $\cZ_s$ and its image in $\Sn$ (via the second projection $\left\lbrace s \right\rbrace \times \Sn \cong \Sn$, which is an isomorphism) by $F_s$. This identification together with Example \ref{exam} gives isomorphisms
\begin{equation}\label{isom}
I_{F_s} \cong \left(\IZ \right)_s \cong \Phi(\mathcal{O}_s). 
\end{equation}

\begin{remark}
	By the definition of $\cZ$ given in \eqref{univ} we have
	\begin{equation*}
		F_s=\left\lbrace [Z]\in \Sn\,|\, s\in\text{supp}(Z) \right\rbrace \subset \Sn.
	\end{equation*}
That is, $F_s$ is the subscheme of $\Sn$ classifying all length $n$ subschemes containing the closed point $s\in S$ in its support. Note that $\dim(F_s)=2n-2$, that is $F_s$ is a subscheme of codimension two in $\Sn$.
\end{remark}

Since $S$ is integral the Hilbert polynomial of $F_s$ does not depend on $s\in S$ by \cite[Proposition 2.1.2]{huy}, call it $p(t)$. We thus have a well defined classifying morphism
\begin{equation*}
 	\varphi: S \rightarrow \Hilb^{p(t)}(\Sn),\,\,\, s\mapsto \left[ F_s\right].
\end{equation*}

\begin{remark}
	Here we can choose any ample line bundle $L\in \Pic(\Sn)$ to define the Hilbert polynomial $p(t)$ of $F_s$, as there is no distinguished ample line bundle on $\Sn$. The choice of a different ample line bundle $\overline{L}$ would give rise to a different Hilbert polynomial $\overline{p}(t)$, but it would not change the proof of the main theorem.
	
	A more conceptual way would be to choose a $n$-very ample line bundle $M$ on $S$, then by \cite{cat} there is a closed embedding
	\begin{equation*}
		\Sn \rightarrow \mathrm{Gr}(n,\mathrm{H}^0(S,M)^{*}),\,\,\, \left[Z \right] \mapsto \mathrm{H}^0(S,M\otimes\mathcal{O}_Z)^{*}.
	\end{equation*}
Composing this morphism with the Plücker embedding of the Grassmannian, we get a closed embedding $\Sn \rightarrow \mathbb{P}^N$ for some $N$ and we can pullback $\mathcal{O}_{\mathbb{P}^N}(1)$ to get an ample line bundle $L$ on $\Sn$. 
\end{remark}

We claim that the morphism $\varphi$ identifies $S$ with a smooth connected component of $\Hilb^{p(t)}(\Sn)$. To see this we have to show that the morphism is injective on closed points and that for every closed point $s\in S$ we have
\begin{equation*}
	\dim(T_{\left[F_s \right] }\Hilb^{p(t)}(\Sn))=2.
\end{equation*}

We start by picking two closed points $s_1\neq s_2\in S$ and note that using equation \eqref{isom} as well as Corollary \ref{ext}, we get:
\begin{equation}\label{Homvan}
	\Hom_{\Sn}(I_{F_{s_1}},I_{F_{s_2}})\cong \Hom_{\Sn}(\Phi(\mathcal{O}_{s_1}),\Phi(\mathcal{O}_{s_2})\cong \Hom_S(\mathcal{O}_{s_1},\mathcal{O}_{s_2})=0.
\end{equation}
If $\left[F_{s_1}\right]= \left[F_{s_2}\right]\in\Hilb^{p(t)}(\Sn)$, then we would have an isomorphism $\mathcal{O}_{F_{s_1}}\cong \mathcal{O}_{F_{s_2}}$ and the exact sequences (for $i=1,2$)
\begin{equation*}
		\begin{tikzcd}
			0 \arrow[r] & I_{F_{s_i}} \arrow[r] & \mathcal{O}_{\Sn} \arrow[r] & \mathcal{O}_{F_{s_i}} \arrow[r] & 0
		\end{tikzcd}
	\end{equation*}
would give rise to a commutative diagram of short exact sequences (with the identity between $\mathcal{O}_{\Sn}$) which shows that there is a nontrivial morphism between $I_{F_{s_1}}$ and $I_{F_{s_2}}$. But this is impossible by \ref{Homvan}. So the classifying morphism $\varphi$ is indeed injective on closed points.

To find $\dim(T_{\left[F_s \right] }\Hilb^{p(t)}(\Sn))$ we remark that
\begin{equation*}
	T_{\left[F_s \right] }\Hilb^{p(t)}(\Sn)\cong\Hom_{\Sn}(I_{F_s},\mathcal{O}_{F_s}),
\end{equation*}
see for example \cite[Proposition 2.2.7]{huy}.

As $q=0$ we have $\Pic^0(S)=0$, but by \cite[Theorem 5.4.]{foga} we also have an isomorphism
\begin{equation*}
	\Pic^0(S)\stackrel{\cong}{\longrightarrow}\Pic^0(\Sn)
\end{equation*}
and thus $\Pic^0(\Sn)=0$. So we can use \cite[Lemma B.5.6.]{kuz} which gives an isomorphism
\begin{equation*}
	\Hom_{\Sn}(I_{F_s},\mathcal{O}_{F_s})\cong \Ext^1_{\Sn}(I_{F_s},I_{F_s}).
\end{equation*}
We find, using again equation \eqref{isom} and Corollary \ref{ext}:
\begin{align*}
	\Ext^1_{\Sn}(I_{F_{s}},I_{F_{s}})&\cong \Ext^1_{\Sn}(\Phi(\mathcal{O}_{s}),\Phi(\mathcal{O}_{s}))\cong \Ext^1_S(\mathcal{O}_{s},\mathcal{O}_{s})\cong T_sS.
\end{align*}
Putting all results together shows $\dim(T_{\left[F_s \right] }\Hilb^{p(t)}(\Sn))=2$ as desired. \qed
\begin{remark}
	Theorem \ref{faith} can be generalized in the case $n=2$ to all smooth projective varieties $X$ having the property $\mathrm{H}^i(X,\mathcal{O}_X)=0$ for $i\geqslant 1$ (that is $\mathcal{O}_X$ is exceptional). By \cite[Theorem A]{bel} the integral functor $\Phi: \Db(X)\rightarrow \Db(X^{[2]})$ with kernel the universal ideal sheaf $\IZ$ is also fully faithful in these cases. Thus our proof of the main result is also valid in these cases.
\end{remark}
\begin{remark}
	In the case of surfaces, the proof of the main result only works for those surfaces with $p_g=q=0$, since for $n\geqslant 2$ the integral functor $\Phi: \Db(S)\rightarrow \Db(\Sn)$ is fully faithful if and only if $p_g=q=0$ by \cite[Theorem A.]{bel2}. But there are similar results for K3 surfaces as well as abelian surfaces, see \cite{rz}. In these cases the integral functor $\Phi$ is a so-called $\mathbb{P}^n$-functor, which again allows to reduce cohomological computations on $\Sn$ to computations on $S$.
\end{remark}


\begin{thebibliography}{10}
	
	\bibitem{bauer2}
	Ingrid Bauer and Fabrizio Catanese.
	\newblock Some new surfaces with {$p_g=q=0$}.
	\newblock In {\em The {F}ano {C}onference}, pages 123--142. Univ. Torino,
	Turin, 2004.
	
	\bibitem{bauer}
	Ingrid Bauer, Fabrizio Catanese, and Roberto Pignatelli.
	\newblock Surfaces of general type with geometric genus zero: a survey.
	\newblock In {\em Complex and differential geometry}, volume~8 of {\em Springer
		Proc. Math.}, pages 1--48. Springer, Heidelberg, 2011.
	
	\bibitem{bel}
	Pieter Belmans, Lie Fu, and Theo Raedschelders.
	\newblock Hilbert squares: derived categories and deformations.
	\newblock {\em Selecta Math. (N.S.)}, 25(3):Paper No. 37, 32, 2019.
	
	\bibitem{bel2}
	Pieter {Belmans} and Andreas {Krug}.
	\newblock Derived categories of (nested) {H}ilbert schemes.
	\newblock {\em arXiv e-prints}, September 2019, to appear in Michigan Mathematical Journal.
	
	\bibitem{bel3}
	Pieter Belmans and Swarnava Mukhopadhyay.
	\newblock Admissible subcategories in derived categories of moduli of vector
	bundles on curves.
	\newblock {\em Adv. Math.}, 351:653--675, 2019.
	
	\bibitem{cat}
	Fabrizio Catanese and Lothar G\"ottsche.
	\newblock {$d$}-{V}ery-ample line bundles and embeddings of {H}ilbert schemes
	of {$0$}-cycles.
	\newblock {\em Manuscripta Math.}, 68(3):337--341, 1990.
	
	\bibitem{dolga}
	Igor Dolgachev.
	\newblock Algebraic surfaces with {$q=p_g=0$}.
	\newblock In {\em Algebraic surfaces}, volume~76 of {\em C.I.M.E. Summer Sch.},
	pages 97--215. Springer, Heidelberg, 2010.
	
	\bibitem{foga}
	John Fogarty.
	\newblock Algebraic families on an algebraic surface. {II}. {T}he {P}icard
	scheme of the punctual {H}ilbert scheme.
	\newblock {\em Amer. J. Math.}, 95:660--687, 1973.
	
	\bibitem{fon}
	Anton Fonarev and Alexander Kuznetsov.
	\newblock Derived categories of curves as components of {F}ano manifolds.
	\newblock {\em J. Lond. Math. Soc. (2)}, 97(1):24--46, 2018.
	
	\bibitem{huy2}
	Daniel Huybrechts.
	\newblock {\em Fourier-{M}ukai transforms in algebraic geometry}.
	\newblock Oxford Mathematical Monographs. The Clarendon Press, Oxford
	University Press, Oxford, 2006.
	
	\bibitem{huy}
	Daniel Huybrechts and Manfred Lehn.
	\newblock {\em The geometry of moduli spaces of sheaves}.
	\newblock Cambridge Mathematical Library. Cambridge University Press,
	Cambridge, second edition, 2010.
	
	\bibitem{krug}
	Andreas Krug and Pawel Sosna.
	\newblock On the derived category of the {H}ilbert scheme of points on an
	{E}nriques surface.
	\newblock {\em Selecta Math. (N.S.)}, 21(4):1339--1360, 2015.
	
	\bibitem{krug1}
	Andreas {Krug} and J{\o}rgen {Vold Rennemo}.
	\newblock {Some ways to reconstruct a sheaf from its tautological image on a
		Hilbert scheme of points}.
	\newblock {\em arXiv:1808.05931}, pages 1--18, 2018.
	\newblock To appear in Math. Nachr.
	
	\bibitem{kuz}
	Alexander Kuznetsov, Yuri Prokhorov, and Constantin Shramov.
	\newblock Hilbert schemes of lines and conics and automorphism groups of {F}ano
	threefolds.
	\newblock {\em Jpn. J. Math.}, 13(1):109--185, 2018.
	
	\bibitem{lang}
	H.~Lange and P.~E. Newstead.
	\newblock On {P}oincar\'{e} bundles of vector bundles on curves.
	\newblock {\em Manuscripta Math.}, 117(2):173--181, 2005.
	
	\bibitem{moon}
	Kyoung-Seog Lee and Han-Bom Moon.
	\newblock Positivity of the {P}oincar\'e bundle on the moduli space of vector
	bundles and its applications.
	\newblock {\em arXiv e-prints}, June 2021.
	
	\bibitem{mum}
	David Mumford.
	\newblock Further pathologies in algebraic geometry.
	\newblock {\em Amer. J. Math.}, 84:642--648, 1962.
	
	\bibitem{nara1}
	M.~S. Narasimhan.
	\newblock Derived categories of moduli spaces of vector bundles on curves.
	\newblock {\em J. Geom. Phys.}, 122:53--58, 2017.
	
	\bibitem{nara2}
	M.~S. Narasimhan.
	\newblock Derived categories of moduli spaces of vector bundles on curves {II}.
	\newblock In {\em Geometry, algebra, number theory, and their information
		technology applications}, volume 251 of {\em Springer Proc. Math. Stat.},
	pages 375--382. Springer, Cham, 2018.
	
	\bibitem{rz}
	Fabian Reede and Ziyu Zhang.
	\newblock Examples of smooth components of moduli spaces of stable sheaves.
	\newblock {\em Manuscripta Math.}, 165(3-4):605--621, 2021.
	
\end{thebibliography}
\end{document}